\documentclass[10pt,a4paper]{article}
\title{\textbf{On paraquaternionic submersions between paraquaternionic K\"ahler manifolds}}
\author{Angelo V. Caldarella}
\date{}

\usepackage{amsfonts}
\usepackage{amsmath,amssymb,amscd,amsthm}
\usepackage{multicol,longtable,tabularx}
\usepackage{enumerate}
\usepackage{latexsym}
\usepackage{amstext,graphics,color}
\usepackage{ifthen}
\usepackage{verbatim}
\usepackage{calc}

\newtheorem{theorem}{Theorem}[section]

\newtheorem{corollary}[theorem]{Corollary}

\newtheorem{proposition}[theorem]{Proposition}
\theoremstyle{definition}
\newtheorem{definition}[theorem]{Definition}

\newtheorem{remark}[theorem]{Remark}

\newcommand{\proofend}{\hfill{\mbox{$\Box$}}\medskip}
\newcommand{\proofbegin}{\noindent{\sc Proof.\ }}
\renewcommand{\proof}{\proofbegin}
\renewcommand{\endproof}{\proofend}

\begin{document}

%
%
%
%
%
%
%
%

\maketitle

\begin{abstract}
In this paper we deal with some properties of a class of semi-Riemannian submersions between manifolds endowed with paraquaternionic structures, proving a result of non-existence of paraquaternionic submersions between paraquaternionic K\"ahler non locally hyper paraK\"ahler manifolds. Then we examine, as an example, the canonical projection of the tangent bundle, endowed with the Sasaki metric, of an almost paraquaternionic Hermitian manifold.
\end{abstract}

\section{Introduction}
The theory of (semi-)Riemannian submersions, as a ``dual'' of that of Riemannian immersions, is a relatively new and vast subject of study, which since its introduction made by O'Neill in \cite{ONArt}, and indipendently by Gray in \cite{G}, has been continously developing, due to its growing importance in the physical framework which try to find a unifying theory between the four fundamental force fields and the gravitational force field. A systematic exposition of the subject may be found in Chapter $9$ of \cite{B}, and in the monographic work \cite{FIP}, to which the reader is referred for more details.

In this paper we introduce a new class of semi-Riemannian submersions, namely the paraquaternionic submersions, which rise in the contest of local differential geometry of manifolds endowed with structures related to the algebra of paraquaternions. We prove, as main result, that the structures on the base space and on the total space of such a class of semi-Riemannian submersion has to be locally hyper paraK\"ahler. After the introductive Section $2$, presenting some preliminar facts on paraquaternionic structures, in Section $3$ we give the definition of paraquaternionic submersion, establishing some basic properties of horizontal and vertical distributions. Section $4$ is devoted to the proof of the main result, and finally, in Section $5$, we present an example of paraquaternionic submersion, namely the canonical projection from the tangent bundle of a paraquaternionic K\"ahler manifold, ending up proving some further properties.

\section{Preliminaries}

We recall some basic data about paraquaternionic manifolds. For a more detailed treatment of the subject, the reader is referred, for example, to \cite{G-RMV-L}, \cite{IMZ}, \cite{IZ} and \cite{V}, and particularly, for the geometry of paracomplex structures, see \cite{AMT} and \cite{CFG}.

\begin{definition}\label{defAPS}
An \emph{almost paraquaternionic structure} on a manifold $M$ is any rank-three subbundle $\sigma $ of the vector
bundle $End(TM)$, such that for each $x\in M$ there exists an open neighbourhood $U$ of $x$ and a
triple $(J_{1},J_{2},J_{3})$ of $(1,1)$--type tensor fields on $U$, satisfying the following conditions:
\begin{enumerate}
  \item[(i)] For any $y\in U$ the triple $((J_{1})_{y},(J_{2})_{y},(J_{3})_{y})$
spans the vector space $\sigma _{y}$
  \item[(ii)] $J_a^2=-\tau_aI$, \quad for any $a\in\{1,2,3\}$;
  \item[(iii)] $J_aJ_b=\tau_cJ_c=-J_bJ_a$, \quad for any cyclic permutation $(a,b,c)$ of $(1,2,3)$;
\end{enumerate}
where $\tau_1=\tau_2=-\tau_3=-1$. Such a kind of triple is said to be a \emph{local basis} (on $U$) for the almost paraquaternionic structure $\sigma$.
The couple $(M,\sigma)$ is said \emph{almost paraquaternionic manifold}.
\end{definition}

When $J_1,J_2,J_3$ are globally defined and satisfy (ii) and (iii) of the previous definition, then the structure $(J_1,J_2,J_3)$ is called \emph{almost hyper paracomplex}.
	
An useful criterion for the construction of a paraquaternionic structure on a smooth manifold is given with the following result.

\begin{proposition}[\cite{C}]
\label{propAPM}
Let $M$ be a manifold, carrying an open covering 
$(U_{h})_{h\in H}$, such that for any $h\in H$ an almost hyper paracomplex
structure $(J_a^{h})_{a=1,2,3}$ is defined on $U_{h}$ satisfying
the following condition: for any $h,k\in H$, with $U_{h}\cap U_{k}\neq
\emptyset $, there exists a smooth map $s_{hk}:U_{h}\cap U_{k}\rightarrow
GL(3,\mathbb{R})$ such that $J_a^{k}|_{U_{h}\cap U_{k}}=(s_{hk})_a^bJ_b^{h}|_{U_{h}\cap U_{k}}$. Then, if $x\in M$, choosing an $h\in H$ such that $x\in U_{h}$, and putting $\sigma _{x}=\mathrm{Span}_{\mathbb{R}}((J_{1}^{h})_{x},(J_{2}^{h})_{x},(J_{3}^{h})_{x})$, we get a three dimensional vector subspace $\sigma_x$ of $End(T_{x}M)$, which does not depend on the open set $U_{h}$ and on the structure $(J_a^{h})_{a=1,2,3}$, such that $\sigma:=\bigcup_{x\in M}\sigma _{x}$ is an almost paraquaternionic structure on $M$.
\end{proposition}

\begin{definition}[\cite{G-RMV-L}, \cite{IZ}]
Let $(M,\sigma)$ be an almost paraquaternionic manifold, and $g\in \mathfrak{T}%
_{2}^{0}(M)$ a metric tensor field. The metric will be said $\sigma$--%
\emph{Hermitian} if the following compatibility condition holds 
\begin{equation}
g(X,J_aY)+g(J_aX,Y)=0,
\label{eq03}
\end{equation}
for any local basis $(J_a)_{a=1,2,3}$ of $\sigma$, defined on an open set $U\subset M$, and any $X,Y\in \Gamma(TU)$. The triple $(M,\sigma,g)$ is said to be an \emph{almost paraquaternionic Hermitian manifold}.

An almost paraquaternionic Hermitian manifold $(M,\sigma,g)$ is said to be \emph{paraquaternionic Hermitian} if there exists on $M$ a torsion-free linear connection which preserves the subbundle $\sigma$. Such a connection is called \emph{paraquaternionic}.

In particular, if on an almost paraquaternionic Hermitian manifold $(M,\sigma,g)$ the Levi-Civita connection is paraquaternionic, then $(M,\sigma,g)$ is said to be \emph{paraquaternionic K\"{a}hler manifold}.
\end{definition}

Obviously, $\dim(M)=4m$, and the metric $g$ is semi-Riemannian, with signature $(2m,2m)$. The following is a well-known characterization for paraquaternionic K\"ahler manifolds.

\begin{proposition}[\cite{G-RMV-L}, \cite{IMZ}, \cite{IZ}]
\label{propPQK02}
Let $(M,\sigma,g)$ be an almost paraquaternionic Hermitian manifold. Then, it is a paraquaternionic K\"{a}hler manifold if and only if for any $x\in M$ there exists an open neighbourhood $U$ of $x$ on which a local basis $(J_a)_{a=1,2,3}$ for $\sigma$ is defined, such that the Levi-Civita connection verifies $\nabla J_a=-\tau_c\omega_c\otimes J_b+\omega_b\otimes J_c$, for any cyclic permutation $(a,b,c)$ of $(1,2,3)$, where $(\omega_a)_{a=1,2,3}$ are 1-forms defined on $U$, and $\tau_{1}=\tau_{2}=-1=-\tau_{3}$.
\end{proposition}

\begin{remark}\label{rm001}
An almost paraquaternionic Hermitian manifold $(M,\sigma,g)$ is paraquaternionic K\"ahler if, and only if, for any local basis $(J_a)_{a=1,2,3}$ for $\sigma$ there exist local 1-forms $(\omega_a)_{a=1,2,3}$, such that the Levi-Civita connection verifies $\nabla J_a=-\tau_c\omega_c\otimes J_b+\omega_b\otimes J_c$, for any cyclic permutation $(a,b,c)$ of $(1,2,3)$,
\end{remark}

\begin{definition}
An almost paraquaternionic Hermitian manifold $(M,\sigma,g)$ is said to be \emph{locally hyper paraK\"{a}hler} if, for any $x\in M$, there exists an open neighbourhood $U$ of $x$, on which a local basis
$(J_{\alpha })_{\alpha=1,2,3}$ of $\sigma$ is defined, such that $\nabla
J_{\alpha }=0$, for each $\alpha\in\{1,2,3\}$.
\end{definition}

Clearly, an almost paraquaternionic Hermitian manifold is locally hyper paraK\"ahler if and only if it is paraquaternionic K\"ahler, with $\omega_a=0$, for any $a\in\{1,2,3\}$.

\begin{definition}
Let $(\bar{M},\bar{\sigma })$ be an almost
paraquaternionic manifold. A submanifold $M$ of $\bar{M}$ is said to be a \emph{$\bar{\sigma}$-invariant
submanifold} of $\bar{M}$ if, for any $x\in M$, and any $\bar{J}\in \bar{\sigma }_{x}$, one has $\bar{J}(T_{x}M)\subset T_{x}M $.
\end{definition}

The above definition is given as the paraquaternionic counterpart of that presented, for the quaternionic setting, in \cite{AM}. If $M$ is a $\bar{\sigma}$-invariant submanifold of an almost paraquaternionic manifold $(\bar{M},\bar{\sigma})$, then it is easy to see that the structure $\bar{\sigma}$ induces, by restriction, an almost paraquaternionic structure $\sigma$ on $M$. Now, let $(M,\sigma)$ be a $\bar{\sigma}$-invariant submanifold of an almost paraquaternionic Hermitian manifold $(\bar{M},\bar{\sigma},\bar{g})$. If we suppose that, for each $x\in M$, the tangent space $T_{x}M$ is a non degenerate subspace of $(T_{x}\bar{M},\bar{g}_{x})$, then necessarily the index of $g_{x}:=\bar{g}_{x}| _{T_{x}M\times T_{x}M}$ does not depend on the point $x$. Therefore, with the only request that each $T_{x}M$ is a non degenerate subspace of $T_{x}\bar{M}$, we get that $(M,g)$ is a semi-Riemannian submanifold of $(\bar{M},\bar{g})$. The induced metric $g$ is obviously $\sigma$--compatible, thus $(M,\sigma,g)$ is an almost paraquaternionic Hermitian manifold, which is said to be an \emph{almost paraquaternionic Hermitian submanifold} of $\bar{M}$.

\begin{proposition}[\cite{C}, \cite{V}]
\label{propPKTGS}Let $(\bar{M},\bar{\sigma },\bar{g})$ be a
paraquaternionic K\"{a}hler manifold, and $(M,\sigma ,g)$ an almost
paraquaternionic Hermitian submanifold of $\bar{M}$. Then $(M,\sigma ,g)$ is also a
paraquaternionic K\"{a}hler manifold, totally geodesic submanifold of $%
\bar{M}$.
\end{proposition}

\section{Paraquaternionic submersions}
Following the analogue definition given in the quaternionic context (cf.\ \cite{IMV01} and \cite{IMV02}), we introduce the following.
\begin{definition}
\label{defPHM}
Let $(M,\sigma )$ and $(M',\sigma')$ be
almost paraquaternionic manifolds. A smooth map $\pi:M\rightarrow M'$ is said to be a \emph{$(\sigma,\sigma')$--{paraholomorphic map}} if it
satisfies the following condition: for any $x\in M$ and for any open neighbourhood $U$ of $x$, on which a local basis $(J_a)_{a=1,2,3}$ for $\sigma$ is defined, there exists an open neighbourhood $U'$ of $\pi(x)$, on which a local basis $(J'_a)_{a=1,2,3}$ for $\sigma'$ is defined, such that $(J'_a)_{\pi(y)}\circ (d\pi)_{y}=(d\pi)_{y}\circ (J_a)_{y}$, for any $y\in U\cap \pi^{-1}(U')$, and any $a\in\{1,2,3\}$.
\end{definition}

It is easy to prove the following basic properties (cf.\ \cite{C}).

\begin{proposition}
\label{propVert}
Let $(M,\sigma ,g)$ and $(M',\sigma ',g')$ be
almost paraquaternionic Hermitian manifolds, and $\pi :M\rightarrow
M'$ a $C^{\infty}$-submersion. If $\pi $ is a $(\sigma
,\sigma ')$--paraholomorphic map, then for each local basis $%
(J_{\alpha })_{\alpha=1,2,3}$ of $\sigma $, we have

\begin{enumerate}
 \item[\emph{(i)}]  The vertical and horizontal distributions $\mathcal{V}$ and $\mathcal{H}$ are both $J_{\alpha
}$--invariant;
 \item[\emph{(ii)}]  The canonical projections $v:\Gamma(TM)\rightarrow\Gamma(\mathcal{V})$ and $h:\Gamma(TM)\rightarrow\Gamma(\mathcal{H})$ commute with each $J_{\alpha }$.
\end{enumerate}
\end{proposition}

As an immediate consequence of the previous result, we have.

\begin{proposition}
Let $(M,\sigma ,g)$ and $(M',\sigma ',g')$ be
almost paraquaternionic Hermitian manifolds, and $\pi :M\rightarrow
M'$ a $C^{\infty}$-submersion. If $\pi $ is a $(\sigma
,\sigma ')$--paraholomorphic map, then each fiber $\pi
^{-1}(x')$, $x'\in M'$, is an almost
paraquaternionic Hermitian submanifold of $(M,\sigma ,g)$.
\end{proposition}

\begin{definition}
Let $(M,\sigma ,g)$ and $(M',\sigma ',g')$ be
almost paraquaternionic Hermitian manifolds, and $\pi :M\rightarrow
M'$ a smooth map. We shall call $\pi $ a \emph{%
paraquaternionic submersion} if it is both a semi--Riemannian submersion and
a $(\sigma ,\sigma ')$--paraholomorphic map.
\end{definition}

\begin{proposition}
\label{propPKM+S}Let $(M,\sigma ,g)$ and $(M',\sigma ^{\prime
},g')$ be almost paraquaternionic Hermitian manifolds, and $\pi
:M\rightarrow M'$ a paraquaternionic submersion. If $M$ is a
paraquaternionic K\"{a}hler manifold, also $M'$ is a
paraquaternionic K\"{a}hler manifold.
\end{proposition}

\proof We claim that, for any $x'\in M'$, we may consider an open neighbourhood $U'$ of $x'$ and an open set $U\subset M$, carrying local bases $(J'_a)_{a=1,2,3}$ and $(J_a)_{a=1,2,3}$ for $\sigma'$ and $\sigma$, respectively, such that $\pi(U)=U'$ and
\begin{equation}\label{Eq002}
J'_a\circ\pi_{*}=\pi_{*}\circ J_a,
\end{equation}
for any $a\in\{1,2,3\}$, where we denote the restriction $(\pi|_{U})_{\sharp}:U\rightarrow U'$ still with $\pi$. Indeed, if $x'\in M'$, choosing $x\in M$ such that $\pi(x)=x'$, and taking an open neighbourhood $\tilde{U}$ of $x$ which carries a local basis $(\tilde{J}_a)_{a=1,2,3}$ for $\sigma$, since $\pi$ is a $(\sigma,\sigma')$-paraholomorphic map, by Definition \ref{defPHM}, there exists an open neighbourhood $\tilde{U}'$ of $x'$ which carries a local basis $(\tilde{J}'_a)_{a=1,2,3}$ for $\sigma'$, such that $(\tilde{J}'_a)_{\pi(y)}\circ (d\pi)_{y}=(d\pi)_{y}\circ (\tilde{J}_a)_{y}$, for any $y\in \tilde{U}\cap \pi^{-1}(\tilde{U}')$ and any $a\in\{1,2,3\}$. Since $\pi$ is a submersion, it is an open map, thus putting $U=\tilde{U}\cap \pi^{-1}(\tilde{U}'\cap\pi(\tilde{U}))$ and $U'=\pi(U)$, and restricting the local bases $(\tilde{J}_a)_{a=1,2,3}$ and $(\tilde{J}'_a)_{a=1,2,3}$ to $U$ and to $U'$, respectively, we get our claim. Since $(M,\sigma,g)$ is a paraquaternionic K\"{a}hler manifold, by Remark \ref{rm001}, 
there exist 1-forms $(\omega_a)_{a=1,2,3}$ on $U$, such that for each cyclic permutation $(a,b,c)$ of $(1,2,3)$: 
\begin{equation}\label{Eq003}
\nabla J_{a}=-\tau _{c}\omega _{c}\otimes J_{b}+\omega _{b}\otimes J_{c},
\end{equation}
where $\tau_{1}=\tau_{2}=-1=-\tau_{3}$. The map $\pi:U\rightarrow U'$ is a semi--Riemannian submersion, 
and for any basic vector fields $X,Y,Z$ on $U$ which are $\pi$-related to $X',Y',Z'\in \Gamma(TU')$ one has that the horizontal component of $\nabla _{X}Y$ is a basic vector field $\pi$-related to $\nabla'_{X'}Y'$ and thus it follows easily that $(\nabla _{X'}'J_{a}')(Y')\circ\pi=\pi _{*}((\nabla _{X}J_{a})(Y))$.
Using (\ref{Eq002}) and (\ref{Eq003}), we have, for any $x\in U$:
\begin{equation}\label{eq31}
(\nabla _{X'}'J_{a}')(Y')_{\pi(x)}=-\tau _{c}\omega _{c}(X)(x)((J'_{b})(Y'))_{\pi(x)}
                                           +\omega _{b}(X)(x)((J_{c}')(Y'))_{\pi(x)},
\end{equation}
from which, by the linear independence of the $J'_{a}$, one gets that the functions $\omega _{a}(X)$, defined on $U$, are constant on
each fiber. This allows us to project each function $\omega _{a}(X)$ on $U'$, that is there exist unique functions $\omega _{a}'(X')$ defined on $U'$, such that $\omega _{a}'(X')\circ \pi =\omega _{a}(X)$. So, we get local 1-forms $(\omega'_{a})_{a=1,2,3}$
on $U'$, simply considering, for any $X'\in \Gamma(TU')$, $\omega'_{\alpha}(X')$ as the unique
projection of $\omega _{a}(X)$ on $U'$, where $X\in \Gamma(TU)$ is the unique basic vector field $\pi$-related to $X'$.
Clearly, by (\ref{eq31}) one has, for each $X',Y'\in\Gamma (TU')$: 
\[
(\nabla _{X'}'J_{a}')(Y')=-\tau_{c}\omega _{c}'(X')J'_{b}Y'+\omega _{b}'(X')J'_{c}Y'
\]
that is, by Proposition \ref{propPQK02}, the manifold $(M',\sigma',g')$ is paraquaternionic
K\"{a}hler. \endproof

Obvious consequences of the previous result are the following.

\begin{corollary}
\label{propHPKM+S}Let $(M,\sigma ,g)$ and $(M',\sigma ^{\prime},g')$ be almost paraquaternionic Hermitian manifolds, and $\pi:M\rightarrow M'$ a paraquaternionic submersion. If $M$ is a locally hyper paraK\"{a}hler manifold, also $M'$ is a locally hyper paraK\"{a}hler manifold.
\end{corollary}

\begin{corollary}
\label{propVertDistr}
Let $(M,\sigma ,g)$ and $(M',\sigma ',g')$ be
paraquaternionic K\"{a}hler manifolds, and $\pi :M\rightarrow M'$ a paraquaternionic submersion. Then the fibers of $\pi $ are totally geodesic paraquaternionic K\"{a}hler submanifolds of $M$.
\end{corollary}

\begin{proposition}\label{corollHIntgr}
Let $(M,\sigma ,g)$ and $(M',\sigma ',g')$ be
paraquaternionic K\"{a}hler manifolds, and $\pi :M\rightarrow M'$ a
paraquaternionic submersion. Then the horizontal distribution is integrable and its integral
manifolds are totally geodesic paraquaternionic K\"{a}hler submanifolds of $M$. 
\end{proposition}

\proof It suffices to prove that the O'Neill tensor $A$ vanishes, and to this purpose we have only to prove that $A=0$
on $\Gamma(\mathcal{H})\times\Gamma(\mathcal{H})$. To this end, let $%
E,F\in \Gamma (\mathcal{H})$. Since $M$ is a paraquaternionic
K\"{a}hler manifold, we may choose a local basis $(J_{a})_{a=1,2,3}$ for $\sigma $ on an open set $U\subset M$, for which there exist 1-forms $(\omega _{a})_{a=1,2,3}$ on $U$, such that, for any $V\in \Gamma (\mathcal{V})$, and any cyclic permutation $(a,b,c)$ of $(1,2,3)$, using the $J_{a}$-invariance of $\mathcal{H}$, one gets $g((\nabla _{E}J_{\alpha })(F),V)=g(-\tau _{g}\omega _{g}(E)J_{b}F+\omega _{b}(E)J_{g}F,V)=0$. Therefore, $(\nabla_{E}J_{\alpha })(F)\in \Gamma (\mathcal{H})$. Then using $v\circ J_{a}=J_{a}\circ v$, by the definition of the O'Neill tensor field $A$, for any $E,F\in \Gamma (\mathcal{H})$, and any $a\in \left\{
1,2,3\right\}$ one gets $(A_{E}J_{a})(F))=v((\nabla _{E}J_{a})(F))=0$, hence $A_{E}(J_{a}F)=J_{a}(A_{E}F)$.
Using the anticommutativity of $A$ on the horizontal distribution, one easily get also, for any $E,F\in \Gamma (\mathcal{H})$, and any $a\in
\left\{ 1,2,3\right\}$, that $A_{J_{a}E}F=J_{a}(A_{E}F)$.
Finally, using (ii) and (iii) of Definition \ref{defAPS} and the previous two relations, one has
\[
A_{E}F=-A_{E}((J_{3})^{2}(F))=-J_{3}(A_{E}(J_{1}J_{2}F))=-J_{3}(A_{J_{2}J_{1}E}F)=(J_{3})^{2}(A_{E}F)=-A_{E}F,
\]
for any $E,F\in \Gamma (\mathcal{H})$, hence $A=0$ on $\Gamma (\mathcal{H})\times \Gamma (\mathcal{H})$, concluding the proof. \endproof 

We recall the following well-known result (see, for example, \cite{BF}).

\begin{proposition}
Let $(M,\tilde{g})$ be a semi--Riemannian manifold endowed with two complementary orthogonal non degenerate integrable distributions $\mathcal{D}$ and $\mathcal{D}^{\bot}$, whose integral manifolds are totally geodesic submanifolds of $M$. Denoting with $\mathcal{F}$ and $\mathcal{F}^{\bot}$ the foliations on $M$ defined, respectively, by $\mathcal{D}$ and $\mathcal{D}^{\bot}$, for any $x^{*}\in M$, there exists an open neighbourhood $U^{*}\subset M$ of $x^{*}$, and two open submanifolds $U$ and $U^{\bot}$ respectively of the leaf of $\mathcal{F}$ and of $\mathcal{F}^{\bot}$ passing through $x^{*}$, such that, denoting with $g$ and $g^{\bot}$ the metrics naturally induced on $U$ and $U^{\bot}$, $(U,\tilde{g}|_{U})$ is the semi--Riemannian product of $(U,g)$ and $(U^{\bot},g^{\bot})$.
\end{proposition}

\begin{remark}
In the hypotheses of Proposition \ref{corollHIntgr}, we have that $M$ is locally a semi--Riemannian product and, as we shall see in the next section, this forces $M$ and $M'$ to be locally hyper paraK\"ahler. To this purpose we need to establish a preliminary property of paraquaternionic K\"ahler manifolds that are also locally semi--Riemannian product, from the point of view of integrable distributions.
\end{remark}

\section{Paraquaternionic K\"ahler structures and semi-Rieman-nian products}

\begin{definition}
Let $(M,\sigma)$ be an almost paraquaternionic manifold, and $\mathcal{D}$ a distribution on $M$. We say that the distribution $\mathcal{D}$ is \emph{$\sigma$--invariant} if, for any $x\in M$ and any $J\in \Sigma_{x}$, one has $J(\mathcal{D}_{x})\subset\mathcal{D}_{x}$.
\end{definition}

It is easy to see that a distribution $\mathcal{D}$ on an almost paraquaternionic manifold $(M,\sigma)$ is $\sigma$--invariant if and only if for any open set $U\subset M$ carrying a local basis $(J_{a})_{a=1,2,3}$ for $\sigma$, one has $(J_{a})_{y}(\mathcal{D}_{y})\subset\mathcal{D}_{y}$, for any $y\in U$. Moreover if the $\sigma$--invariant distribution $\cal{D}$ is integrable, then each integral submanifold is an almost paraquaternionic submanifold of $(M,\sigma)$. In particular, if $(M,\sigma,g)$ is a paraquaternionic K\"ahler manifold, and $\cal{D}$ is a non degenerate, integrable, $\sigma$--invariant distribution, then each integral manifold of $\sigma$ is a totally geodesic paraquaternionic K\"ahler submanifold of $(M,\sigma,g)$.

Let us suppose now that an almost paraquaternionic Hermitian manifold $(M,\sigma,g)$ is endowed with two non degenerate complementary orthogonal distributions $\cal{D}$ and $\cal{D}'$; then, using (\ref{eq03}), one easily has that the $\sigma$--invariance of one of the distributions implies the $\sigma$--invariance of the other.

Now, we are interested in the study of a paraquaternionic K\"ahler manifold endowed with two non degenerate complementary orthogonal, $\sigma$--invariant, integrable distributions $\cal{D}$ and $\cal{D}'$. Indeed, each integral manifold of both $\cal{D}$ and $\cal{D}'$ is a totally geodesic paraquaternionic K\"ahler submanifold of $(M,\sigma,g)$, then $M$ is a locally semi--Riemannian product. To this aim, we need some basic facts about the product of two almost paraquaternionic structures, which we are going to report, following partially the exposition made in \cite{V}.

Let $(M_{1},\sigma_{1},g_{1})$ and $(M_{2},\sigma_{2},g_{2})$ be two almost paraquaternionic Hermitian manifolds, and put $M=M_{1}\times M_{2}$. Consider the semi--Riemannian product metric $g=g_{1}\times g_{2}$ on $M$. Given two open sets $U_{1}\subset M_{1}$ and $U_{2}\subset M_{2}$, on which local bases $(J^{(1)}_{a})_{a=1,2,3}$ and $(J^{(2)}_{a})_{a=1,2,3}$ for $\sigma_{1}$ and $\sigma_{2}$, respectively, are defined, then on $U:=U_{1}\times U_{2}\subset M$ we define, for any $a\in\{1,2,3\}$, a $(1,1)$--type tensor field $J_{a}$ by putting, for any $X\in\Gamma(TM)$: $J_{a}X:=J^{(1)}_{a}P_{1}X+J^{(2)}_{a}P_{2}X$, where $P_1$ and $P_2$ are the canonical projections of $\Gamma(TM)$ onto $\Gamma(TM_1)$ and $\Gamma(TM_2)$, respectively. It is known (cf.\ \cite{V}) that the triple $(J_{a})_{a=1,2,3}$ is an almost hyper paracomplex structure on the open set $U$, which we shall call \emph{of product type}. Supposing $(M_{1},\sigma_{1},g_{1})$ and $(M_{2},\sigma_{2},g_{2})$ be paraquaternionic K\"ahler manifolds, then there exist 1-forms $(\omega^{(1)}_{a})_{a=1,2,3}$ on $U_{1}$, and 1-forms $(\omega^{(2)}_{a})_{a=1,2,3}$ on $U_{2}$, such that, if $\nabla^{(1)}$ and $\nabla^{(2)}$ are the Levi-Civita connections on $M_{1}$ and $M_{2}$, one has $\nabla^{(\lambda)}J_{a}^{(\lambda)}=-\tau_{c}\omega^{(\lambda)}_{c}\otimes J^{(\lambda)}_{b}+\omega^{(\lambda)}_{b}\otimes J^{(\lambda)}_{c}$, for any $\lambda\in\{1,2\}$ and any cyclic permutation $(a,b,c)$ of $(1,2,3)$, with $\tau_{1}=\tau_{2}=-1=-\tau_{3}$. Furthermore, defining 1-forms $(\omega_{a})_{a=1,2,3}$ on $U$  by putting $\omega_{a}(X):=\omega^{(1)}_{a}(P_{1}X)+\omega^{(2)}_{a}(P_{2}X)$, then the following identities hold, for any $X,Y\in\Gamma(TU)$ and any cyclic permutation $(a,b,c)$ of $(1,2,3)$ (cf. \cite[Lemma 3.3]{V}):
\begin{equation}
\label{eq45}
2(\nabla_{X}J_{a})Y=-\tau_{c}\omega_{c}(X)J_{b}Y+\omega_{b}(X)J_{c}Y-\tau_{c}\omega_{c}(FX)J_{b}(FY)+\omega_{b}(FX)J_{c}(FY),
\end{equation}
where $\tau_{1}=\tau_{2}=-1=-\tau_{3}$ and $F:=P_{1}-P_{2}$.

Now we are able to prove the following result.

\begin{theorem}
Let $(M,\sigma,g)$ be a paraquaternionic K\"ahler manifold, endowed with two non degenerate, complementary orthogonal, integrable, $\sigma$--invariant distributions $\mathcal{D}_{1}$ and $\mathcal{D}_{2}$. Then $(M,\sigma,g)$ is a locally hyper paraK\"ahler manifold.
\end{theorem}
\proof Being $M$ locally a semi--Riemannian product, if $x\in M$, then we may consider an open neighbourhood $U\subset M$ of $x$, for which there exist two open submanifolds $U_{1}$ and $U_{2}$ of the leaves $L^{(1)}_{x}$ and $L^{(2)}_{x}$, respectively, of the foliations defined by $\mathcal{D}_{1}$ and $\mathcal{D}_{2}$, passing through the point $x$, such that, denoting with $g_{1}$ and $g_{2}$, respectively, the metrics induced from $g$ on $L^{(1)}_{x}$ and $L^{(2)}_{x}$, the manifold $(U,g)$ is isometric to the semi--Riemannian product of the manifolds $(U_{1},g_{1})$ and $(U_{2},g_{2})$.
Since $U$ is an open submanifold of $(M,\sigma,g)$, it has a paraquaternionic K\"ahler structure, which we shall denote with $\sigma'$, naturally induced, by restriction, from $\sigma$, so that $(U,\sigma',g)$ is a paraquaternionic K\"ahler manifold. Each local basis $(J'_{a})_{a=1,2,3}$ of $\sigma'$, defined on an open set $U'\subset U$, is obtained by taking a local basis $(J_{a})_{a=1,2,3}$ of $\sigma$ on an open set $\bar{U}'\subset M$, such that $U'=\bar{U}'\cap U$, and putting $J'_{a}:=J_{a}|_{U'}$, $a\in\{1,2,3\}$. Being $U'$ an open set of $M$, it is clear that $(J'_{a})_{a=1,2,3}$ is not only a local basis for $\sigma'$, but also a local basis for $\sigma$. Therefore, we may say that the structure $\sigma'$ is spanned by the family of all the local bases $(J'_{a})_{a=1,2,3}$ for $\sigma$, defined on open sets $W\subset M$, provided that $W\subset U$.

Since the distributions $\mathcal{D}_1$ and $\mathcal{D}_2$ are both $\sigma$-invariant, then the leaves $L^{(1)}_{x}$ and $L^{(2)}_{x}$ are both $\sigma$-invariant submanifolds of $M$, therefore they carry paraquaternionic K\"ahler structures, which we shall denote respectively with $\sigma_{1}$ and $\sigma_{2}$, naturally induced from the structure $\sigma$. By Proposition \ref{propPKTGS}, we have that $(L^{(1)}_{x},\sigma_{1},g_{1})$ and $(L^{(2)}_{x},\sigma_{2},g_{2})$ are totally geodesic paraquaternionic K\"ahler submanifolds of $(M,\sigma,g)$. As before, the open submanifolds $(U_{1},g_{1})$ and $(U_{2},g_{2})$ inherit, in a natural way, paraquaternionic K\"ahler structures, $\sigma'_{1}$ and $\sigma'_{2}$, from $\Sigma_{1}$ and $\Sigma_{2}$, respectively.

Now, let us consider a local basis $(J'_a)_{a=1,2,3}$ for $\sigma'$, defined on an open neighbourhood $U'\subset U$ of $x$. It is also a local basis for $\sigma$, and we may consider local bases $(J'_{a}{}^{(1)})_{a=1,2,3}$ and $(J'_{a}{}^{(2)})_{a=1,2,3}$ for $\sigma'_1$ and $\sigma'_2$, respectively, naturally induced from $(J'_a)_{a=1,2,3}$ as follows. Since $U\cong U_1\times U_2$ and $U'$ is an open subset of $U$, up to a suitable restriction, we may take open sets $U'_1\subset U_1$ and $U'_2\subset U_2$, such that $U'\cong U'_1\times U'_2$.
Furthermore, being $U_1\subset L^{(1)}_x$ and $U_2\subset L^{(2)}_x$, we have $U'_1=U'\cap L^{(1)}_x$ and $U'_2=U'\cap L^{(2)}_x$, Thus, if we set $J'_{a}{}^{(1)}:=(J'_a)|_{U'_1}$ and $J'_{a}{}^{(2)}:=(J'_a)|_{U'_2}$, for any $a\in\{1,2,3\}$, we get local bases $(J'_{a}{}^{(1)})_{a=1,2,3}$ for $\Sigma_1$, and $(J'_{a}{}^{(2)})_{a=1,2,3}$ for $\Sigma_2$, defined on the open sets $U'_1\subset L^{(1)}_x$ and $U'_2\subset L^{(2)}_x$, respectively. Being $U'_1\subset U_1$ and $U'_2\subset U_2$, we can say that $(J'_{a}{}^{(1)})_{a=1,2,3}$ and $(J'_{a}{}^{(2)})_{a=1,2,3}$ are also local bases for $\sigma'_1$ and $\sigma'_2$, respectively.

Beside the almost hyper paracomplex structure $(J'_{a})_{a=1,2,3}$ fixed at the beginning on $U'$, we may define on it another almost hyper paracomplex structure. Namely, since $U'\cong U'_{1}\times U'_{2}$, we consider the almost hyper paracomplex structure $(J''_{a})_{a=1,2,3}$ of product type, obtained from $(J'_{a}{}^{(1)})_{a=1,2,3}$ and $(J'_{a}{}^{(2)})_{a=1,2,3}$. By the definition of $J'_{a}{}^{(1)}$ and $J'_{a}{}^{(2)}$, denoting with $P_1:\Gamma(TU')\rightarrow \Gamma(TU'_1)$ and $P_2:\Gamma(TU')\rightarrow \Gamma(TU'_2)$ the canonical projections, one has $J''_{a}X:=J'_{a}{}^{(1)}(P_{1}X)+J'_{a}{}^{(2)}(P_{2}X)=J'_{a}(P_{1}X)+J'_{a}(P_{2}X)=J'_{a}X$, for any $X\in\Gamma(TU')$, and any $a\in\{1,2,3\}$. We find, in this way, that the local basis $(J'_{a})_{a=1,2,3}$, chosen at the beginning, is of product type.

Being $(M,\sigma,g)$ a paraquaternionic K\"ahler manifold, and $(J'_{a})_{a=1,2,3}$ a local basis for $\sigma$, there exist 1--forms $(\omega_{a})_{a=1,2,3}$ on $U'$, such that
\begin{equation}
\label{eq43}
(\nabla_{X}J'_{a})Y=-\tau_{c}\omega_{c}(X)J'_{b}Y+\omega_{b}(X)J'_{c}Y
\end{equation}
for any $X,Y\in\Gamma(TU')$ and any cyclic permutation $(a,b,c)$ of $(1,2,3)$, with $\tau_{1}=\tau_{2}=-1=-\tau_{3}$. At the same time, $(U_{1},\sigma'_{1},g_{1})$ and $(U_{2},\sigma'_{2},g_{2})$ are open submanifolds of the paraquaternionic K\"ahler manifolds $(L^{(1)}_{x},\sigma_{1},g_{1})$ and $(L^{(2)}_{x},\sigma_{2},g_{2})$, respectively, and so they are themselves paraquaternionic K\"ahler manifolds. Being $(J'_{a}{}^{(1)})_{a=1,2,3}$ and $(J'_{a}{}^{(2)})_{a=1,2,3}$ local bases for $\sigma'_{1}$ and $\sigma'_{2}$, respectively, there exist 1--forms $(\omega_{a}^{(1)})_{a=1,2,3}$ defined on $U'_{1}$, and $(\omega_{a}^{(2)})_{a=1,2,3}$ defined on $U'_{2}$, such that $\nabla^{(\lambda)}J'_{a}{}^{(\lambda)}=-\tau_{c}\omega^{(\lambda)}_{c}\otimes J'_{b}{}^{(\lambda)}+\omega^{(\lambda)}_{b}\otimes J'_{c}{}^{(\lambda)}$,
for any $\lambda\in\{1,2\}$, and any cyclic permutation $(a,b,c)$ of $(1,2,3)$. Since $(J'_{a})_{a=1,2,3}$ is a structure of product type, if we define on $U'$ 1--forms $(\bar{\omega}_{a})_{a=1,2,3}$ by putting, for any $a\in\{1,2,3\}$ and any $X\in\Gamma(TU')$, $\bar{\omega}_{a}(X):=\omega_{a}^{(1)}(P_{1}X)+\omega_{a}^{(2)}(P_{2}X)$, then, by (\ref{eq45}), we have
\begin{equation} \label{eq46}
2(\nabla_{X}J'_{a})Y=-\tau_{c}\bar{\omega}_{c}(X)J'_{b}Y+\bar{\omega}_{b}(X)J'_{c}Y-\tau_{c}\bar{\omega}_{c}(FX)J'_{b}(FY)+\bar{\omega}_{b}(FX)J'_{c}(FY),
\end{equation}
for any $X,Y\in\Gamma(TU')$ and any cyclic permutation $(a,b,c)$ of $(1,2,3)$, with $\tau_{1}=\tau_{2}=-1=-\tau_{3}$ and $F:=P_{1}-P_{2}$. Since the local bases $(J'_{a}{}^{(1)})_{a=1,2,3}$ and $(J'_{a}{}^{(2)})_{a=1,2,3}$ are defined by restriction of the local basis $(J'_{a})_{a=1,2,3}$ then also the 1--forms $(\omega_{a}^{(1)})_{a=1,2,3}$ and $(\omega_{a}^{(2)})_{a=1,2,3}$ are obtained by restriction of the 1--forms $(\omega_{a})_{a=1,2,3}$. Therefore, one easily obtains that $\bar{\omega}_{a}=\omega_{a}$, for any $a\in\{1,2,3\}$. So, (\ref{eq46}) may be rewritten as follows
\begin{equation}\label{eq47}
2(\nabla_{X}J'_{a})Y=-\tau_{c}\omega_{c}(X)J'_{b}Y+\omega_{b}(X)J'_{c}Y-\tau_{c}\omega_{c}(FX)J'_{b}(FY)+\omega_{b}(FX)J'_{c}(FY)
\end{equation}
for any $X,Y\in\Gamma(TU')$ and any cyclic permutation $(a,b,c)$ of $(1,2,3)$. Comparing (\ref{eq43}) and (\ref{eq47}), and decomposing $X$ and $Y$ along $U'_1$ and $U'_2$, we get
\[
\tau_{c}\omega_{c}(P_{1}X)J'_{b}(P_{2}Y)+\tau_{c}\omega_{c}(P_{2}X)J'_{b}(P_{1}Y)-\omega_{b}(P_{1}X)J'_{c}(P_{2}Y)-\omega_{b}(P_{2}X)J'_{c}(P_{1}Y)=0
\]
for any $X,Y\in\Gamma(TU')$, and any $b,c\in\{1,2,3\}$, with $b\neq c$. Using the $\sigma$--invariance of the distributions $\mathcal{D}_{1}$ and $\mathcal{D}_{2}$, and comparing the components along $\mathcal{D}_{1}$ and $\mathcal{D}_{2}$, we get
\[
\left\{\begin{array}{l}
\tau_{c}\omega_{c}(P_{1}X)J'_{b}(P_{2}Y)-\omega_{b}(P_{1}X)J'_{c}(P_{2}Y)=0 \\
\tau_{c}\omega_{c}(P_{2}X)J'_{b}(P_{1}Y)-\omega_{b}(P_{2}X)J'_{c}(P_{1}Y)=0
\end{array}\right.,
\]
for any $X,Y\in\Gamma(TU')$ and any $b,c\in\{1,2,3\}$, with $b\neq c$. Using the definition of $J'_{a}{}^{(1)}$ and $J'_{a}{}^{(2)}$, we have
\[
\left\{\begin{array}{l}
\tau_{c}\omega_{c}(P_{1}X)J'_{b}{}^{(2)}-\omega_{b}(P_{1}X)J'_{c}{}^{(2)}=0 \\
\tau_{c}\omega_{c}(P_{2}X)J'_{b}{}^{(1)}-\omega_{b}(P_{2}X)J'_{c}{}^{(1)}=0
\end{array}\right.,
\]
for any $X\in\Gamma(TU')$ and any $b,c\in\{1,2,3\}$, with $b\neq c$, from which, since $(J'_{a}{}^{(1)})_{a=1,2,3}$ and $(J'_{a}{}^{(2)})_{a=1,2,3}$ are almost hyper paracomplex structures on $U'_{1}$ and $U'_{2}$, respectively, it follows that $\omega_{a}(P_{1}X)|_{U'_{1}}=0$ and $\omega_{a}(P_{2}X)|_{U'_{2}}=0$. This implies $\omega_{a}=0$, for any $a\in\{1,2,3\}$, finally getting $\nabla J'_{a}=0$. \endproof

\begin{proposition}
Let $(M,\sigma,g)$ and $(M',\sigma',g')$ be paraquaternionic K\"ahler manifolds, and $\pi: M\rightarrow M'$ a paraquaternionic submersion. Then, both $(M,\sigma,g)$ and $(M',\sigma',g')$ are locally hyper paraK\"ahler manifolds.
\end{proposition}
\proof Since $\pi$ is a paraquaternionic submersion, by Proposition \ref{propVert}, Corollary \ref{propVertDistr} and Corollary \ref{corollHIntgr}, the vertical and the horizontal distributions are two non degenerate, complementary orthogonal, integrable, $\sigma$--invariant distributions on $M$. Thus, the manifold $(M,\sigma,g)$ is locally hyper paraK\"ahler, and by Corollary \ref{propHPKM+S}, also the manifold $(M',\sigma',g')$ is locally hyper paraK\"ahler. \endproof

Let us now give a brief and simple reformulation of the above results in terms of (almost) product structures. It is well known (\cite{AM}, \cite{BF}, \cite{CFG}) that the assignment of two complementary distributions $\mathcal{D}$ and $\mathcal{D}'$ on a manifold $M$ is equivalent to the assignment of an \emph{almost product structure} on it, that is a $(1,1)$--type tensor field $F\neq\pm I$, satisfying $F^2=I$, with $\mathcal{D}=T^+M$ and $\mathcal{D}'=T^-M$, where $T^+M$ and $T^-M$ are the eigendistributions with respect the eigenvalues $\pm1$ of $F$. The integrability of both distributions is equivalent to the vanishing of the Nijenhuis tensor field $N_F$ related to the structure $F$ (see \cite{Y}), and in this case the tensor field $F$ is called a \emph{product structure}, or a \emph{locally product structure}.

An \emph{indefinite Riemannian almost product structure} (cf.\ \cite{ES}) on $M$ is a pair $(g,F)$, where $g$ is a semi--Riemannian metric, and $F$ is an almost product structure satisfying $g(FX,FY)=g(X,Y)$, for any $X,Y\in\Gamma(TM)$. The triple $(M,g,F)$ is called an \emph{indefinite Riemannian almost product manifold}. We shall drop the term ''almost'' from the previous denominations, when $N_{F}=0$. In this case, the eigendistributions $T^{\pm}M$ are complementary orthogonal distributions.

It is easy to prove the following result.

\begin{proposition}
\label{prop004}
Let $(M,g,F)$ be an indefinite Riemannian product manifold. The integral manifolds of the eigendistributions $T^{+}M$ and $T^{-}M$ are totally geodesic submanifolds of $M$ if and only if $\nabla F=0$, being $\nabla$ the Levi-Civita connection on $M$.
\end{proposition}

Let us introduce the following definition.

\begin{definition}
Let $(M,\sigma)$ be an almost paraquaternionic manifold. An almost product structure $F$ on $M$ is said to be \emph{$\sigma$--invariant} if for any open set $U$ endowed with a local basis $(J_a)_{a=1,2,3}$ for $\sigma$, one has $J_a\circ F|_U=F|_U\circ J_a$, for any $a\in\{1,2,3\}$.
\end{definition}

\begin{proposition}
Let $(M,\sigma,g)$ be a paraquaternionic K\"ahler manifold, and $F$ a $\sigma$--invariant almost product structure on $M$, such that $(M,g,F)$ is an indefinite Riemannian almost product manifold. The following three conditions are equivalent:
\begin{enumerate}
 \item[\emph{(i)}]$\nabla F=0$:
 \item[\emph{(ii)}]$N_F=0$;
 \item[\emph{(iii)}]$(\nabla_{J_aX}F)(Y)-(\nabla_XF)(J_aY)=0$;
\end{enumerate}
where, in (iii), $(J_a)_{a=1,2,3}$ is a local basis on an open set $U$, $a\in\{1,2,3\}$ and $X,Y\in\Gamma(TU)$. Furthermore, if any of the above conditions is verified, then $(M,\sigma,g)$ is a locally hyper paraK\"ahler manifold.
\end{proposition}
\proof Since the equivalence between (i) and (ii), as well as the last statement, are immediate, we only prove the equivalence between (i) and (iii). If (i) holds, then (iii) clearly holds. Choose an open set $U$ endowed with a local basis $(J_a)_{a=1,2,3}$ for $\sigma$. Since $M$ is a paraquaternionic K\"ahler manifold, using the $\sigma$-invariance of $F$, one gets
\begin{equation}\label{eq:001}
F(\nabla_XJ_a)(Y)=(\nabla_XJ_a)(FY)	
\end{equation}
for any $a\in\{1,2,3\}$, and any $X,Y\in\Gamma(TU)$. Since $\nabla_X(F\circ J_a)(Y)=(\nabla_XF)(J_aY)+F(\nabla_XJ_a)(Y)$ and $\nabla_X(J_a\circ F)(Y)=(\nabla_XJ_a)(FY)+J_a(\nabla_XF)(Y)$, using again the $\sigma$-invariance of $F$ and (\ref{eq:001}), one has $(\nabla_XF)(J_aY)=J_a(\nabla_XF)(Y)$, for any $a\in\{1,2,3\}$, and any $X,Y\in\Gamma(TU)$. If we suppose that (iii) holds, then for any $a\in\{1,2,3\}$, and any $X,Y\in\Gamma(TU)$, we get $(\nabla_XF)(J_aY)=J_a(\nabla_XF)(Y)=(\nabla_{J_aX}F)(Y)$, thus
{\setlength\arraycolsep{2pt}
\begin{eqnarray*}
(\nabla_{X}F)(Y)&=&-(\nabla_XF)(J_3^2Y)=-J_3(\nabla_XF)(J_1J_2Y)\\
                &=&-J_3(\nabla_{J_2J_1X}F)(Y)=J_3^{2}(\nabla_XF)(Y)=-(\nabla_XF)(Y),
\end{eqnarray*}}
from which $\nabla F=0$ follows. \endproof

\section{An example of paraquaternionic submersion}
Let $M$ be an $m$-dimensional manifold and $(TM,\pi,M)$ its tangent bundle. We recall the following basic properties of the vertical and horizontal lifts of vector fields, following \cite{D}, \cite{I}, \cite{K} and \cite{YS}. If $\nabla$ is a linear
connection on $M$, then we have, for any $X,Y\in \Gamma (TM)$, and any $\xi \in TM$: $\left[ X^{v},Y^{v}\right] _{\xi }=0$, $\left[ X^{h},Y^{v}\right] _{\xi }=(\nabla _{X}Y)_{\xi }^{v}$ and $\left[ X^{h},Y^{h}\right] _{\xi }=-R(X,Y,Z)_{\xi }^{v}+\left[X,Y\right] _{\xi }^{h}$, being $Z\in \Gamma (TM)$ arbitrarily chosen, such that $Z_{\pi(\xi )}=\xi $, and $R$ the $(1,3)$--type curvature tensor field of
the connection $\nabla $. Let $g$ be a semi--Riemannian metric on $M$, and $\nabla$ its the Levi-Civita connection. The \emph{Sasaki metric} on $TM$ is the semi--Riemannian metric $G$ determined by $G(X^{v},Y^{v})=g(X,Y)\circ \pi$, $G(X^{h},Y^{h})=g(X,Y)\circ\pi$ and $G(X^{v},Y^{h})=0$, for any $X,Y\in \Gamma (TM)$. If $(p,q)$ is the signature of the metric $g$, then $(2p,2q)$ is the signature of the Sasaki metric $G$. Then, denoting with $\widetilde{\nabla}$ the Levi-Civita connection of the Sasaki metric, for any $X,Y\in \Gamma (TM)$, and any 
$\xi \in TM$, we have
{\setlength\arraycolsep{2pt}
\begin{eqnarray}
&&(\widetilde{\nabla }_{X^{v}}Y^{v})_{\xi }=0, \qquad (\widetilde{\nabla }_{X^{h}}Y^{h})_{\xi }=-\frac{1}{2}R(X,Y,Z)_{\xi}^{v}+(\nabla _{X}Y)_{\xi }^{h}, \label{eq01}\\
&&(\widetilde{\nabla }_{X^{h}}Y^{v})_{\xi }=(\nabla _{X}Y)_{\xi }^{v}+\frac{1}{2}R(Z,Y,X)_{\xi }^{h}, \qquad (\widetilde{\nabla }_{X^{v}}Y^{h})_{\xi }=\frac{1}{2}R(Z,X,Y)_{\xi}^{h},\label{eq02}
\end{eqnarray}}being $Z\in \Gamma (TM)$ arbitrarily chosen, such that $Z_{\pi(\xi )}=\xi $. Finally, we have

\begin{proposition}[\cite{FIP}]
\label{propFIP}Let $(M,g)$ be a semi--Riemannian manifold. Then the
canonical projection $\pi :(TM,G)\rightarrow (M,g)$ is a semi--Riemannian
submersion with totally geodesic fibers; moreover, the horizontal
distribution $\mathcal{H}TM$ is integrable if and only if the metric $g$ is
flat, and in this case its integral manifolds are totally geodesic submanifolds of $TM$.
\end{proposition}

Let $(M,\sigma ,g)$ be an almost paraquaternionic Hermitian manifold, and endow the tangent bundle $TM$ with the Sasaki metric. Let us choose an open covering $(U_{i})_{i\in I}$ of $M$, such that, for any $i\in I$, a local basis $(J^{i}_{a})_{a=1,2,3}$ of hyper paracomplex structures for $\sigma$ is defined on $U_i$. Let us define on $\pi^{-1}(U_{i})$, for any $i\in I$, the $(1,1)$--type tensor fields $(\tilde{J}^{i}_{1},\tilde{J}^{i}_{2},\tilde{J}^{i}_{3})$, putting, for any $X\in \Gamma (TU)$, and any $a \in \left\{1,2,3\right\}$
\begin{equation}
\tilde{J}^{i}_{a }(X^{v}):=(J^{i}_{a }(X))^{v}\quad \quad \text{and}%
\quad \quad \tilde{J}^{i}_{a }(X^{h}):=(J^{i}_{a }(X))^{h}.  \label{eq08}
\end{equation}

\begin{proposition}
Let $(M,\sigma,g)$ be an almost paraquaternionic Hermitian manifold, and endow the tangent bundle $TM$ with the Sasaki metric. The almost paraquaternionic structure $\sigma$ induces, via $\pi$, an almost paraquaternionic structure $\tilde{\sigma}$ on $TM$, such that $(TM,\tilde{\sigma},G)$ is an almost paraquaternionic Hermitian manifold, and $\pi: TM\rightarrow M$ is a paraquaternionic submersion. 
\end{proposition}
\proof
As seen above, we get an open covering $(\pi^{-1}(U_{i}))_{i\in I}$ of $TM$, such that, for any $i\in I$, an almost hyper paracomplex structure $(\tilde{J}^{i}_{a})_{a=1,2,3}$ is defined on the open set $\pi^{-1}(U_{i})$. If $i,j\in I$ are such that $U_{i}\cap U_{j}\neq\emptyset$, there exists a smooth map $s_{ij}: U_{i}\cap U_{j}\rightarrow GL(3,\mathbb{R})$, such that $J^{i}_{a}|_{U_{i}\cap U_{j}}=(s_{ij})^{b}_{a}J^{j}_{b}|_{U_{i}\cap U_{j}}$, for any $a\in\{1,2,3\}$. If $X\in\Gamma(T(U_{i}\cap U_{j}))$, using (\ref{eq08}) and the above equality, we have $\tilde{J}^{i}_{a}(X^{v})=((s_{ij})^{b}_{a}\circ\pi)\tilde{J}^{j}_{b}(X^{v})$, for any $a\in\{1,2,3\}$. Analogously, for $X^{h}$. Putting $\tilde{s}_{ij}:=s_{ij}\circ\pi$, we get a smooth map $\tilde{s}_{ij}: \pi^{-1}(U_{i})\cap\pi^{-1}(U_{j})\rightarrow GL(3,\mathbb{R})$, such that $\tilde{J}^{i}_{a}|_{\pi^{-1}(U_{i})\cap\pi^{-1}(U_{j})}=(\tilde{s}_{ij})^{b}_{a}\tilde{J}^{j}_{b}|_{\pi^{-1}(U_{i})\cap\pi^{-1}(U_{j})}$, for any $a\in\{1,2,3\}$. By Proposition \ref{propAPM}, we conclude that the open covering $(\pi^{-1}(U_{i}))_{i\in I}$ of $TM$, together with the family of almost hyper paracomplex structures $((\tilde{J}^{i}_{a})_{a=1,2,3})_{i\in I}$, defines an almost paraquaternionic structure $\tilde{\sigma}$ on $TM$.

Let us now prove that the Sasaki metric $G$ is $\tilde{\sigma}$--Hermitian. To this aim, let us consider an open set $U\subset M$ on which a local basis $(J_{a})_{a=1,2,3}$ for $\sigma$ is defined, and let $(\tilde{J}_{a})_{a=1,2,3}$ be the induced local basis for $\tilde{\sigma}$ on $\pi^{-1}(U)$. If $X,Y\in\Gamma(TU)$, we have
\[
G(\tilde{J}_{a}(X^{v}),Y^{v})=g(J_{a}X,Y)\circ\pi=-g(X,J_{a}Y)\circ\pi=-G(X^{v},\tilde{J}_{a}(Y^{v}))
\]
and analogously for the horizontal lifts. Clearly, we have also $G(\tilde{J}_{a}(X^{h}),Y^{v})=0=-G(X^{h},\tilde{J}_{a}(Y^{v}))$,  and it follows that $G$ is a $\tilde{\sigma}$--Hermitian metric, that is $(TM,\tilde{\sigma},G)$ is an almost paraquaternionic Hermitian manifold.

Finally, let us prove that $\pi: TM\rightarrow M$ is a paraquaternionic submersion. Being $\pi$ a semi--Riemannian submersion, one has only to prove that it is a $(\tilde{\sigma},\sigma)$--paraholomorphic map. To this end, let us consider a point $\xi\in TM$. By the definition of the structure $\tilde{\sigma}$ induced on $(TM,G)$ from $\sigma$, there exists an open neighbourhood $U$ of $\pi(\xi)$, on which a local basis $(J_{a})_{a=1,2,3}$ for $\sigma$ is defined, such that on $\pi^{-1}(U)$ the local basis $(\tilde{J}_{a})_{a=1,2,3}$ for $\tilde{\sigma}$ is defined by (\ref{eq08}). Any local basis of $\tilde{\sigma}$ is defined in such way. Clearly, one has, for any $\xi'\in\pi^{-1}(U)$, any $X\in\Gamma(TU)$ and any $a \in \left\{ 1,2,3\right\}$
\[
(d\pi )_{\xi' }((\tilde{J}_{a })_{\xi' }(X_{\xi' }^{v}))=0=(J_{a })_{\pi (\xi' )}((d\pi )_{\xi'}(X_{\xi'}^{v})) 
\]
and
\[
(d\pi )_{\xi' }((\tilde{J}_{a })_{\xi' }(X_{\xi' }^{h}))=(J_{a }X)_{\pi (\xi' )}=(J_{a })_{\pi (\xi')}((d\pi )_{\xi'}(X_{\xi' }^{h}))
\]
This proves that $(J_{a})_{\pi(\xi')}\circ(d\pi)_{\xi'}=(d\pi)_{\xi'}\circ(\tilde{J}_{a})_{\xi'}$, for any $a \in \left\{ 1,2,3\right\}$, and any $\xi'\in\tilde{U}$, and, by Definition \ref{defPHM}, $\pi$ is a $(\tilde{\sigma},\sigma)$--paraholomorphic map.\endproof

\begin{proposition}
Let $(M,\sigma ,g)$ be an almost paraquaternionic Hermitian
manifold. Then, $(TM,\tilde{\sigma},G)$ is a paraquaternionic
K\"{a}hler manifold if and only if $(M,\sigma ,g)$ is a paraquaternionic
K\"{a}hler manifold, with $g$ flat; moreover, also the metric $G$ is flat.
\end{proposition}

\proof Assume that $(TM,\tilde{\sigma},G)$ is a paraquaternionic
K\"{a}hler manifold. Since $\pi$ is a paraquaternionic submersion, then Proposition \ref{propPKM+S},
Corollary \ref{corollHIntgr} and Proposition \ref{propFIP} imply that $(M,\sigma,g)$ is a paraquaternionic K\"{a}hler manifold, with $g$ flat.

Viceversa, let $\xi_{0}\in TM$, and put $\pi(\xi_{0})=x$. Let us consider an open neighbourhood $U\subset M$ of $x$, on which a local basis $(J_{a})_{a=1,2,3}$ for $\sigma$ is defined, and then let us take the local basis $(\tilde{J}_{a})_{a=1,2,3}$ for $\tilde{\sigma}$ on $\pi^{-1}(U)$, induced from $(J_{a})_{a=1,2,3}$. Being $M$ a paraquaternionic K\"{a}hler manifold, there exist 1--forms $(\omega _{a })_{a=1,2,3}$ defined on $U$, such that
\begin{equation}\label{eq07}
\nabla J_{a }=-\tau _{c}\omega _{c}\otimes J_{b}+\omega _{b}\otimes J_{c},
\end{equation}
for any cyclic permutation $(a,b,c)$ of $(1,2,3)$. Using this identity, (\ref{eq01}) and (\ref{eq02}), with
simple calculations one has the following relations, for any $a \in
\left\{ 1,2,3\right\} $%
\begin{equation}
\begin{array}{l}
(\tilde{\nabla }_{X^{v}}\tilde{J}_{a })(Y^{v})_{\xi }=0, \qquad (\tilde{\nabla }_{X^{v}}\tilde{J}_{a })(Y^{h})_{\xi }=\frac{1}{2%
}\left\{ R(Z,X,J_{a }Y)_{\xi }^{h}-J_{a }R(Z,X,Y)_{\xi
}^{h}\right\} ,
\\ 
\\ 
(\tilde{\nabla }_{X^{h}}\tilde{J}_{a })(Y^{h})_{\xi }=-\frac{1}{2%
}\left\{ R(X,J_{a }Y,Z)_{\xi }^{v}-J_{a }R(X,Y,Z)_{\xi
}^{v}\right\} +(\nabla _{X}J_{a })(Y)_{\xi }^{h}, \\ 
\\ 
(\tilde{\nabla }_{X^{h}}\tilde{J}_{a })(Y^{v})_{\xi }=(\nabla
_{X}J_{a })(Y)_{\xi }^{v}+\frac{1}{2}\left\{ R(Z,J_{a }Y,X)_{\xi
}^{h}-J_{a }R(Z,Y,X)_{\xi }^{h}\right\} , 
\end{array}
\label{eq06}
\end{equation}
for any $X,Y\in \Gamma (TM)$, any $\xi \in TM$, and $Z\in \Gamma (TM)$
arbitrarily chosen such that $Z_{\pi (\xi )}=\xi $. Since $g$ is flat, using (\ref{eq08}) and (\ref{eq07}), from (\ref{eq06}) it follows: 
\begin{equation}
\begin{array}{l}
(\tilde{\nabla }_{X^{v}}\tilde{J}_{a })(Y^{v})_{\xi }=0, \qquad (\tilde{\nabla }_{X^{v}}\tilde{J}_{a })(Y^{h})_{\xi }=0,\\\\
(\tilde{\nabla }_{X^{h}}\tilde{J}_{a })(Y^{h})_{\xi }=-\tau
_{c}\omega _{c}(X)_{\pi (\xi )}(\tilde{J}_{b})(Y^{h})_{\xi }+\omega _{b}(X)_{\pi (\xi )}(\tilde{J}_{c})(Y^{h})_{\xi },\\ 
\\ 
(\tilde{\nabla }_{X^{h}}\tilde{J}_{a })(Y^{v})_{\xi }=-\tau
_{c}\omega _{c}(X)_{\pi (\xi )}(\tilde{J}_{b})(Y^{v})_{\xi }+\omega _{b}(X)_{\pi (\xi )}(\tilde{J}_{c})(Y^{v})_{\xi },
\end{array}
\label{eq09}
\end{equation}
Finally, putting, $\tilde{\omega }_{a }:=\pi ^{*}(\omega _{a })$, for any $a \in \left\{ 1,2,3\right\}$, we get 1--forms $(\tilde{\omega }_{a })_{a=1,2,3}$ on $\pi^{-1}(U)$, such that, using (\ref{eq09}), $\tilde{\nabla }\tilde{J}_{a }=-\tau _{c}\tilde{\omega }_{c}\otimes \tilde{J}_{b}+\tilde{\omega }_{b}\otimes \tilde{J}_{c}$, for any cyclic permutation $(a,b,c)$ of $(1,2,3)$. Thus $(TM,\tilde{\sigma },G)$ is a paraquaternionic K\"{a}hler manifold.

Now, let us suppose that one of the two equivalent conditions is satisfied; if $\tilde{R}$ is the curvature of $\tilde{\nabla }$,
then recovering the formulas for the curvature given in \cite{K}, we get that if $R=0$, then immediately $\tilde{R}=0$. \endproof 

\begin{theorem}
Let $(M,\sigma ,g)$ be an almost paraquaternionic Hermitian
manifold. Then, $(TM,\tilde{\sigma},G)$ is locally hyper paraK\"{a}hler if and only if $(M,\sigma,g)$ is locally hyper paraK\"{a}hler, with $g$ flat.
\end{theorem}
\proof Supposing $(TM,\tilde{\sigma},G)$ locally hyper paraK\"{a}hler, since $\pi:TM\rightarrow M$ is a paraquaternionic submersion, by Corollary \ref{propHPKM+S}, we conclude that $(M,\sigma,g)$ is locally hyper paraK\"{a}hler. Further, since $(TM,\tilde{\sigma},G)$ is also a paraquaternionic K\"ahler manifold, we have that $g$ is flat.
Conversely, suppose that $(M,\sigma,g)$ is locally hyper paraK\"ahler, with $g$ flat. Take $\xi_{0}\in TM$, and put $\pi(\xi_{0})=x$. Let us consider an open neighbourhood $U\subset M$ of $x$, on which a local basis $(J_{a})_{a=1,2,3}$ for $\sigma$ is defined, and then let us take the local basis $(\tilde{J}_{a})_{a=1,2,3}$ for $\tilde{\sigma}$ on $\pi^{-1}(U)$, induced from $(J_{a})_{a=1,2,3}$. We may suppose that $\nabla J_{a}=0$, for any $a\in\{1,2,3\}$. Using this condition and the flatness of $g$, by (\ref{eq06}) we get $\tilde{\nabla}\tilde{J}_{a}=0$, for any $a\in\{1,2,3\}$. Thus $(TM,\tilde{\sigma},G)$ is locally hyper paraK\"{a}hler. \endproof

\vspace{10pt}
\noindent \textbf{Acknowledgements}. The author would like to express his gratitude to Prof.\ A.\ M.\ Pastore for her valuable comments and suggestions on the subject, and to Prof.\ S. Ianu\c{s} for many discussions about the topic of the present paper, occurred during his stay at the University of Bari and during the author's visit at the University of Bucharest.

\noindent\small{Author's address}
\vspace{2pt}\\
  \small{Department of Mathematics, University of Bari} \\
  \small{Via E. Orabona 4,}\\
  \small{I-70125 Bari (Italy)}
\vspace{2pt}\\
  \small{\texttt{caldarella@dm.uniba.it}}\\

\end{document}